\documentclass{mcom-l}

\usepackage{lipsum}
\usepackage{amsfonts}
\usepackage{bbold}
\usepackage{graphicx}
\usepackage{epstopdf}
\usepackage{algorithmic}
\usepackage{amsmath}    
\usepackage{amssymb}
\usepackage{graphicx}   
\usepackage{verbatim}   
\usepackage{color}      
\usepackage{hyperref}   
\usepackage{enumerate}
\usepackage{mathrsfs}
\usepackage{empheq}

\newcommand{\del}{\partial}
\renewcommand{\theta}{\vartheta}
\renewcommand{\phi}{\varphi}

\newcommand{\veccc}[3]{\left ( \begin{array}{c}#1\\#2\\#3\\ \end{array}\right )}

\newcommand{\dd}{\text{d}}

\newcommand{\grad}{\text{grad\,}}
\renewcommand{\div}{\text{div\,}}
\renewcommand{\vec}{\mathbf}

\newcommand{\vecomega}{\boldsymbol\omega}

\newcommand{\const}{\text{const}}

\newcommand{\ii}{\mathbb{i}}
\newcommand{\new}[1]{#1}

\newtheorem{theorem}{Theorem}[section]

\theoremstyle{definition}
\newtheorem{definition}[theorem]{Definition}
\newtheorem{example}[theorem]{Example}

\newtheorem{corollary}{Corollary}

\theoremstyle{remark}

\numberwithin{equation}{section}

\begin{document}

\title{Stationarity preserving schemes for multi-dimensional linear systems}

\author{Wasilij Barsukow}

\address{W\"urzburg University, Institute for Mathematics, Emil-Fischer-Stra\ss e 40, 97074 W\"urzburg, Germany}

    \email{w.barsukow@mathematik.uni-wuerzburg.de}

    \thanks{The author acknowledges inspiring discussions with Philip L. Roe and support from the German National Academic Foundation.}

\date{January 0, 0000 and, in revised form, January 00, 0000.}


\begin{abstract}
There is a qualitative difference between one-dimensional and multi-dimensional solutions to the Euler equations: new features that arise are vorticity and a nontrivial incompressible (low Mach number) limit. They present challenges to finite volume methods. 
It seems that an important step in this direction is to first study the new features for the multi-dimensional acoustic equations.
There exists an analogue of the low Mach number limit for this system and its vorticity is stationary.

It is shown that a scheme that possesses a stationary discrete vorticity (vorticity preserving) also has stationary states that are discretizations of all the analytic stationary states. This property is termed stationarity preserving. Both these features are not generically fulfilled by finite volume schemes; in this paper a condition is derived that determines whether a scheme is stationarity preserving (or, equivalently, vorticity preserving) on a Cartesian grid. 

Additionally, this paper also uncovers a previously unknown connection to schemes that comply with the low Mach number limit. Truly multi-dimensional schemes are found to arise naturally and it is shown that a multi-dimensional discrete divergence previously discussed in the literature is the only possible stationarity-preserving one (in a certain class).
\end{abstract}

\maketitle

\section{Introduction}

There is a rich body of work related to finite volume methods for the Euler equations in one spatial dimension (e.g. \cite{roe81}, \cite{harten83} as well as the references in \cite{toro09}, \cite{leveque02} and many others). In multiple spatial dimensions the Euler equations read
\begin{align}
 \del_t \rho + \div(\rho \vec v ) &= 0\\
 \del_t (\rho \vec v) + \div (\rho \vec v \otimes \vec v) + \nabla p &= 0\\
 \del_t e + \div (\vec v (e+p)) &= 0
\end{align}
with
\begin{align}
 e &= \frac{p}{\gamma-1} + \frac12 \rho |\vec v|^2
\end{align}
Here $\rho, \rho \vec v, e$ are the mass, momentum and energy volume densities. $p$ is the pressure and the equation of state of the ideal gas contains the constant $\gamma > 1$.

Solving the Euler equations in multiple spatial dimensions numerically is made difficult by a number of challenges. Finite volume methods can suffer from a bad resolution of vortices and show excessive diffusion in the incompressible (low Mach number) limit (see e.g. \cite{morton01}, \cite{turkel87}, \cite{dellacherie10}, \cite{barsukow16}). Both are truly multi-dimensional features absent in one-dimensional flow. 

In \cite{jeltsch06}, \cite{lerat07}, \cite{mishra09preprint} examples of schemes have been presented that fulfill particular conditions on the evolution of discrete counterparts to the vorticity $\vecomega = \nabla \times \vec v$. This paper is a contribution to the question, how a discrete vorticity is to be chosen and what its desired numerical evolution should be. In order to do this the focus lies on a particular set of equations.

A path that circumvents solving the Euler equations directly has been used in \cite{godlewski13}, \cite{chalons13}, \cite{roe17} and others: they separate the advective operator from the acoustic one. Indeed, already at PDE level the acoustic operator is very different from the advective operator, and it might be advantageous to treat them differently at discrete level. Therefore first a numerical discretization of the acoustic part of the Euler equations is needed. This has led \cite{morton01}, \cite{mishra09preprint}, \cite{amadori2015}, \cite{dellacherierieper10}, \cite{franck17} and others to studies of the linearized form of the acoustic operator. Already for this linear system many features of the behaviour of numerical methods for the nonlinear problem can be seen. The linearized system possesses a stationary vorticity which makes analysis particularly easy. The suggested schemes thus keep a discrete vorticity operator stationary. 

So far only particular examples of vorticity preserving schemes for the acoustic system have been derived (e.g. in \cite{lukacova00}, \cite{morton01}, \cite{jeltsch06}, \cite{mishra09preprint}, \cite{lung14}), and the authors followed very different approaches while deriving the schemes. In this paper the focus lies on finding properties which are common to all vorticity preserving schemes for the acoustic equations. What the examples that were suggested so far have in common is that the vorticity-preserving schemes were all constructed using multi-dimensional stencils, i.e. stencils that (on Cartesian grids in two dimensions, say) involve all the 8 neighbors of a cell rather than just 4. This raises the question of whether this is necessary or whether simpler stencils also allow for vorticity preservation.  Also, the suggested schemes had a very similar appearance, and this paper gives reasons for that.

Additionally, one might want to check, for a given scheme for the acoustic equations, whether it possesses any stationary discrete vorticity operator. In \cite{morton01} the Appendix deals with the preservation of a {specific} discrete vorticity for a certain family of schemes, and with the vorticity production rate in case of non-preservation. However this analysis, as well as the treatment of the acoustic equations in \cite{lerat07}, \emph{assumes} this vorticity discretization. There might still exist some other discrete vorticity that is exactly preserved. One thus needs to be able to check the existence of \emph{any} stationary vorticity discretization. This paper presents solutions to these questions and thus makes possible to explore the entire class of vorticity preserving schemes for the acoustic equations.

Thus, the aim of this paper is to derive general conditions on schemes for the linearized acoustic system to preserve a discrete vorticity. Vorticity preservation turns out to be equivalent to another property of numerical schemes, which is interesting in its own right -- \emph{stationarity preservation}. A numerical scheme is called {stationarity preserving} if it has a discrete counterpart to all the analytic stationary states. Here, stationary states are considered in multiple spatial dimensions, which makes them very diverse. As will be seen later for instance, the acoustic system keeps every divergence-free velocity field stationary. The quest for a discrete counterpart thus makes appear a discrete divergence operator. 

The stationarity-preserving property itself turns out to be tightly linked to the limit of low Mach number; in fact for the acoustic equations the low Mach number limit can be rewritten as the large time limit. This highlights the importance of a correct treatment of stationary states (or vorticity) at discrete level and brings together two different lines of thought: vorticity preserving schemes and schemes that comply with limit of low Mach number.

As is shown in section \ref{sec:example-dimsplit} neither of these properties is a generic feature of finite volume schemes. They achieve stability by introducing diffusion, which is found to destroy all but the very trivial stationary states and to diffuse away vorticity. These schemes therefore only keep a very poor subset of states stationary, and diffuse away all the others. On the other hand, schemes that are stationarity preserving keep stationary a discrete analogue of all the stationary states.
 
{Numerical methods that take special care of particular stationary states are often called \emph{well-balanced}. Although this is not generally equivalent to stationarity preservation, in one spatial dimension relations between well-balancing and the behaviour of these schemes in particular limits have been studied in \cite{gosse15, gosse12, gosse11, gosse02}}.

The paper is structured as follows. In section \ref{sec:fourier-cont} the stationary states of a general linear hyperbolic system of PDEs are discussed by means of the Fourier transform, and it is shown how these concepts immediately translate to the discrete situation. With this preparation in section \ref{sec:statpreserv} a definition of a \emph{stationarity-preserving scheme} is given.

Section \ref{sec:time} discusses why it suffices to leave the time continuous and deal with semi-discrete schemes, and how the results immediately transfer to the fully discrete situation. Here stability of the scheme when subject to explicit time integration is briefly considered.

In section \ref{sec:acoustic} the concepts are applied to the acoustic equations obtained as a linearization of Euler equations (system wave equation). Having a discrete vorticity that is exactly preserved during the numerical time evolution (vorticity preserving scheme) is shown in section \ref{sec:vorticity} to be equivalent to a stationarity preserving scheme. The connection between stationarity preserving and low Mach number schemes for these equations is studied in section \ref{sec:vorticity} and \ref{sec:lowmach}.

Finally, section \ref{sec:consistentdiffusion} shows a way how to construct stationarity-preserving schemes. The procedure is exemplified for first order schemes that solve the acoustic equations. In particular it is shown how, given any discrete divergence, a numerical diffusion is to be constructed that -- loosely speaking -- complies with the discrete divergence being stationary. 
Such a construction is performed explicitly starting from a very simple discrete divergence operator that uses central derivatives. It is shown that the only way to obtain a stationarity preserving (or, equivalently, a vorticity preserving) scheme is to extend this stencil to a truly multi-dimensional one, i.e. it is not possible to find a dimensionally split scheme with these properties. Moreover the corresponding velocity stencils in their most compact forms turn out to be unique if they shall be symmetric. This explains why schemes previously found in literature, although they were derived using very different ideas, have the same essential ingredient: the {\sl factorizable} stencil in \cite{sidilkover02} and the vorticity discretizations obtained in \cite{morton01}, \cite{jeltsch06}, \cite{mishra09preprint}, \cite{lukacova00}. It is just the unique stationarity-preserving discretization for symmetric first-order schemes.

In the conclusions in section \ref{sec:conclusions} an outlook to the nonlinear Euler equations is given.

\section{Stationary states}

First, in this section stationary states of both the linear hyperbolic systems and their discretizations are discussed. Nontrivial stationary states are interesting for numerics because they turn out to be the key to understanding many more properties, like the low Mach number limit and vorticity preservation, which are subjects of later sections. The main result for the discrete situation is Theorem \ref{thm:detcondition}.

\subsection{Continuous case}
\label{sec:fourier-cont}

This section deals with stationary states for the general hyperbolic linear $n \times  n$ system in $d$ spatial dimensions ($\vec J$ being a $d$-dimensional vector of matrices $(J_x, J_y, \ldots)$)
\begin{align}
 \del_t q + \vec J \cdot \nabla q &= 0 \label{eq:schemegeneral} &  q &: \mathbb R^+_0 \times \mathbb R^d \to \mathbb R^n
\end{align}
Although this analysis is inspired by a particular example of such system, namely the acoustic equations discussed in section \ref{sec:acoustic}, stationary states (both at continuous and discrete level) can be fruitfully studied for more general problems. This section is thus devoted to stationary states of general linear systems of the form \eqref{eq:schemegeneral}.

\begin{definition}[Trivial stationary state]
 A stationary state of \eqref{eq:schemegeneral} is called \emph{trivial}, if no component of $q$ can be chosen freely as a function of $\vec x$.
\end{definition}

All the stationary states are characterized by $\vec J \cdot \nabla q$. For example, data that satisfy $\nabla q = 0$ remain stationary for all times. This is a trivial stationary state, because every component of $q$ is a constant. However, for particular choices of $\vec J$ there exist stationary states with more freedom. Recall first that, given any $0 \neq \vec k \in \mathbb R^d$, hyperbolicity of the system \eqref{eq:schemegeneral} guarantees diagonalizability of $\vec J \cdot \vec k$ with real eigenvalues. An eigenvalue zero precisely corresponds to the existence of a richer set of stationary states:
\begin{theorem} \label{thm:stationarycontinuous}
 If $\det (\vec k \cdot \vec J)$ {vanishes for all $\vec k \in \mathbb R^d$}, then there exist non-trivial stationary states of \eqref{eq:schemegeneral}.
\end{theorem}
\begin{proof}
First the vanishing determinant is related to stationary states.

Consider the Fourier transform of \eqref{eq:schemegeneral} by inserting ($\ii = \sqrt{-1}$)
\begin{align}
 q(t, \vec x) = \hat q \exp(- \ii \omega t + \ii \vec k \cdot \vec x)
\end{align}
to obtain the eigenvalue problem $\omega \hat q = \vec J \cdot \vec k \hat q$. 
Every Fourier mode evolves in time as $\exp(- \ii \omega t)$ with the corresponding eigenvalue $\omega$. Initial data that remain stationary have $\omega = 0$; their Fourier transform therefore is the eigenvector $\hat q_0$ of $\vec J \cdot \vec k$ which corresponds to an eigenvalue zero: $\hat q_0 \in \ker (\vec J \cdot \vec k)$. 

Thus it remains to be shown that among these stationary states are nontrivial ones. Consider initial data decomposed in their Fourier modes. For stationarity every Fourier mode $\hat q$ needs to be parallel to $\hat q_0 \in \mathbb R^n$. This is just one condition for $n$ unknowns. Therefore for systems of PDEs, i.e. for $n > 1$, for every $\vec k$ at least one of the components of $\hat q$ can be chosen freely. This implies that also one of the components of $q$ can be chosen freely as a function of space.
\end{proof}

\textit{Remark}:  Scalar problems ($n=1$) only have trivial stationary states.

Trivial stationary states typically do not pose as much challenges to the numerical schemes as non-trivial ones. The remainder of this paper therefore focuses its attention onto systems for which $\det (\vec J \cdot \vec k) = 0 \,\,\forall \vec k$, i.e. those that possess non-trivial stationary states. This is the case for the acoustic equations to be considered later (Corollary \ref{cor:stationarypdelevelacoustic} in section \ref{sec:acoustic}).

\begin{theorem} \label{thm:continuousvorticity}
 The existence of nontrivial stationary states for the system \eqref{eq:schemegeneral} is equivalent to the existence of a constant of motion (one for each zero eigenvalue), i.e. a function $\Omega q : \mathbb R^+_0 \times \mathbb R^d \to \mathbb R$, linear in the solution $q$, which does not evolve in time for any initial data:
 \begin{align}
  \del_t (\Omega q) = 0
 \end{align}
\end{theorem}
\begin{proof}
  Taking the Fourier transform in spatial coordinates only, i.e. inserting 
  \begin{align}
    q(t, \vec x) = \hat q(t) \exp(\ii \vec k \cdot \vec x)
  \end{align}
  yields $\del_t \hat q(t) = - \ii (\vec k \cdot \vec J) \hat q(t)$. Assume that $\det (\vec k \cdot \vec J) = 0$ for all $\vec k$, such that there exists a left eigenvector $\hat \Omega$ that belongs to the eigenvalue zero:
  {
  \begin{align}
   \hat \Omega (\vec k \cdot \vec J) &= 0 &&\Rightarrow  & \del_t (\hat \Omega \hat q) &= 0
  \end{align}
  This eigenvector in general depends on $\vec k$.} Take $\Omega q$ to be the inverse Fourier transform of $\hat \Omega \hat q$. {Note that any factor of $\vec k$ that multiplies $\hat q$ now becomes a derivative that acts on $q$. $\Omega$ is thus in general a differential operator. Then $\del_t (\hat \Omega \hat q) = 0$ for any $\hat q$ and any $\vec k$ implies $\del_t (\Omega q) = 0$.}
\end{proof}

{This constant of motion is referred to as an \emph{involution} in \cite{dafermos86}, where it is emphasized that $\Omega q$ remains zero, if it is zero initially. Here the focus rather lies on the fact that $\Omega q$ is stationary, whatever the initial data. One way to extend these ideas to the discrete case has been presented in \cite{kemm13}, but it is very different from the approach followed below (Section \ref{sec:statpreserv}).}

An example of such a constant of motion for the acoustic equations is given in Corollary \ref{cor:vorticitypdelevelacoustic} of section \ref{sec:acoustic}.

\subsection{Discrete case}
\label{sec:statpreserv}

The study of stationary states by means of the Fourier transform in the proof of Theorem \ref{thm:stationarycontinuous} has a natural equivalent in the discrete sense.

Assume Equation \eqref{eq:schemegeneral} to be solved numerically on a rectangular $d$-dimensional grid.
\begin{definition}[Notation]

 \begin{enumerate}[i)]
  \item The computational grid consists of cells of constant spacing $\Delta x_m$ into the $m$-th spatial direction, $m = 1, \ldots, d$.
  \item Each cell has a unique index which is a $d$-dimensional vector $I \in \mathbb Z^d$ of integer numbers, with components $I_m$, $m = 1, \ldots, d$.
  \item $q_I$ is the value of $q$ in cell $I$.
  \item Given $\vec k \in \mathbb R^d$, its components are denoted by $k_m$, $m = 1, \ldots, d$.
 \end{enumerate}
\end{definition}
Note that in this paper indices never denote derivatives.

The Fourier ansatz now reads
\begin{align}
 q_I = \hat q \exp\left(- \ii \omega t + \ii \sum_{m=1}^d I_m k_m \Delta x_m \right) 
\end{align}

\begin{example} In 2-d one has $I = (i, j)$ and $\Delta x_1 =: \Delta x$, $\Delta x_2 =: \Delta y$, and therefore
\begin{align}
 q_{ij} = \hat q \exp\left(\ii [-\omega t + i k_x \Delta x + j k_y \Delta y] \right)
\end{align}
\hfill$\lhd$\end{example}

\begin{definition}[Translation factor]
The shift by one cell is conveyed by the \emph{translation factor} $t_m(k_m) := \exp(\ii k_m \Delta x_m)$.
\end{definition}
{The dependence on $\vec k$ will be suppressed in the notation.}

This allows to write
\begin{align}
 q_I 
 = \hat q \exp\left(- \ii \omega t\right) \prod_{m=1}^d t_m^{I_m} \label{eq:translationfactor}
\end{align}

Any linear finite difference formula at cell $I$ can be written as
\begin{align}
 \sum_{S \in [-N, N]^d \subset \mathbb Z^d } \alpha_S q_{I + S} \label{eq:generalstencil}
\end{align}

\begin{example}On a two-dimensional grid the central difference in $x$-direction is
\begin{align}
 \frac1{2\Delta x} \left( q_{i+1, j} - q_{i-1, j} \right )
\end{align}
such that $N=1$ and $\alpha_{1,0} = \frac1{2\Delta x}$, $\alpha_{-1,0} = -\frac1{2\Delta x}$ and all other $\alpha$ vanish. {Using \eqref{eq:translationfactor} this becomes
\begin{align}
 \exp(-\ii \omega t ) t_x^i t_y^j \cdot \frac1{2\Delta x} \left( t_x - \frac{1}{t_x} \right ) \hat q
\end{align}
This exemplifies how $t_x$ acts as a shift by one cell to the right, and $t_x^{-1}$ as a shift by one cell to the left.} \hfill$\lhd$\end{example}

{%
\begin{example}The following two-dimensional finite difference formula approximates $\del_x q$:
\begin{align}
 \frac1{8\Delta x} \Big( (q_{i+1, j+1} + 2q_{i+1,j} + q_{i+1,j-1}) - (q_{i-1, j+1} + 2 q_{i-1, j} + q_{i-1, j-1}) \Big )
\end{align}
Its Fourier transform is, up to the usual prefactor $\exp(-\ii \omega t ) t_x^i t_y^j$,
\begin{align}
  \frac1{8\Delta x} \left( t_x - \frac{1}{t_x} \right ) \left( t_y + 2 + \frac{1}{t_y}  \right ) \hat q = \frac{(t_x+1)(t_x-1)}{2t_x \Delta x} \cdot \frac{(t_y+1)^2}{4 t_y} \hat q
\end{align}
\hfill$\lhd$\end{example}
}

In general of course, as the object of study are systems of equations, $q$ is a vector, and every $\alpha_S$ is an $n \times n$ matrix.

Applying the Fourier transform to the finite difference scheme gives rise to a discrete analogue of the condition $\det (\vec k \cdot \vec J)$ found in Theorem \ref{thm:stationarycontinuous}. The role of $\vec k \cdot \vec J$ is played by the \emph{evolution matrix} $\mathcal E$:

\begin{definition}[Evolution matrix]
 The \emph{evolution matrix} associated to the finite difference scheme \eqref{eq:generalstencil} is the matrix
 \begin{align}
\mathcal E(\vec k) \,\,\,\,  = -\!\!\!\!\sum_{S \in [-N, N]^d  } \!\!\!\!\! \ii \alpha_S  \prod_{m=1}^d t_m^{S_m} 
\end{align}
\end{definition}

{Recall that a semi-discrete scheme is called consistent, if for smooth solutions, it converges to the PDE as $\Delta x_m \to 0$, $m = 1, \ldots, d$.}

\begin{definition}[Stationarity preservation]
 A consistent linear scheme for \eqref{eq:schemegeneral} is called \emph{stationarity preserving} if it possesses nontrivial stationary states that discretize all the analytic stationary states of \eqref{eq:schemegeneral}. 
\end{definition} 

{At the continuous level nontrivial stationary states appear if $\vec J \cdot \vec k$ has a vanishing eigenvalue (see proof of Theorem \ref{thm:stationarycontinuous}). Diagonalizability of $\vec J \cdot \vec k$ is guaranteed by hyperbolicity of the system of PDEs. In the discrete case the evolution matrix $\mathcal E$ replaces $\vec J \cdot \vec k$. By analogy in the following it is assumed that $\mathcal E$ is diagonalizable as well (but not necessarily with real eigenvalues).}

\begin{theorem}[Stationarity preservation] \label{thm:detcondition}
 A consistent linear scheme 
 \begin{align}
  \del_t q_I \,\,\,\,  + \!\!\!\! \sum_{S \in [-N, N]^d \subset \mathbb Z^d } \alpha_S q_{I + S} &= 0 \label{eq:semidiscrete}
\end{align}
 with evolution matrix $\mathcal E(\vec k)$ is \emph{stationarity preserving} if for all $\vec k$
 \begin{align}
  \dim \ker \mathcal E(\vec k) = \dim \ker (\vec J \cdot \vec k)
 \end{align}
 A necessary condition is that the determinant $\det \mathcal E(\vec k)$ of $\mathcal E(\vec k)$ vanishes for all $\vec k$. The numerical stationary states are discretizations of the stationary states of the PDE. 
\end{theorem}
\begin{proof}
The finite difference scheme \eqref{eq:generalstencil}, inserting \eqref{eq:translationfactor}, has the Fourier transform
\begin{align}
  \exp\left(- \ii \omega t\right) \left( \prod_{m=1}^d t_m^{I_m} \right ) \!\!\!\! \sum_{S \in [-N, N]^d  } \!\!\!\!\!\alpha_S \left(  \prod_{m=1}^d t_m^{S_m} \right ) \hat q \label{eq:fouriertrafostencilcomplete}
\end{align}

The semidiscrete scheme, assumed to be a consistent discretization of \eqref{eq:schemegeneral},
\begin{align}
  \del_t q_I \,\,\,\,  + \!\!\!\! \sum_{S \in [-N, N]^d \subset \mathbb Z^d } \alpha_S q_{I + S} &= 0 
\end{align}
after the Fourier transform is taken leads to {$\omega \hat q = \mathcal E(\vec k) \hat q $}. {Since $\mathcal E(\vec k)$ is assumed to be diagonalizable,} from here the argument is exactly the same as in the proof of Theorem \ref{thm:stationarycontinuous}, and the existence of nontrivial stationary states is characterized by $\det \mathcal E(\vec k) = 0 \, \forall \vec k$. They are a discretization of the analytical stationarity condition $\vec J \cdot \nabla q = 0$ by consistency. {Indeed, for a consistent scheme, the evolution matrix $\mathcal E(\vec k)$ converges to $\vec J \cdot \vec k$ as $\Delta x_m \to 0$, $m = 1, \ldots, d$. As both the evolution matrix (for all $\Delta x_m$) and $\vec J \cdot \vec k$ are diagonalizable by assumption, the corresponding eigenvectors converge to the eigenvectors of $\vec J \cdot \vec k$.}
\end{proof}

{
Consider a vector $\hat q_0$ parallel to the eigenvector of $\mathcal E$ with eigenvalue zero: $\mathcal E \hat q_0 = 0$. Thus $\del_t \hat q_0 = 0$. In physical space $q_0$ is the Fourier mode $\hat q_0(t) \exp(\ii k_x \Delta x i + \ii k_y \Delta y j + \ldots)$, which thus is kept stationary by the numerical scheme. As stationarity is guaranteed for all $\vec k$, also any linear combinations of such modes with different $\vec k$ remain stationary.} The numerical stationary states, i.e. those states that are kept exactly stationary by the numerics are given by the eigenvector of $\mathcal E$ corresponding to eigenvalue zero. 

Analogously to Theorem \ref{thm:continuousvorticity}, if $\mathcal E$ is diagonalizable, then together with the right eigenvector (characterizing the numerical stationary states) a left eigenvector exists which yields a numerical constant of motion / involution.

\begin{theorem} \label{thm:discretevorticitygeneral}
 Any stationarity preserving scheme with diagonalizable evolution matrix $\mathcal E$ gives rise to a numerical constant of motion, whose Fourier transform is given by the left eigenvector belonging to eigenvalue zero of $\mathcal E$.
\end{theorem}
Replacing $\vec J \cdot \vec k$ by $\mathcal E(\vec k)$ in Theorem \ref{thm:continuousvorticity} proves the assertion.

\subsection{{Time integration and stability}}
\label{sec:time}

So far, the scheme under consideration has left the time continuous. Adopting the viewpoint of the method of lines, however, the results carry over directly to any numerical time integrator that is able to solve the equation $\del_t q = 0$ exactly. This, generally, is a primary requirement upon any such method, and is fulfilled by the forward Euler or any Runge-Kutta scheme. 

Having identified initial data that are left exactly stationary by the numerics, an important question is what happens to other data. Obviously, for stability all Fourier modes that are not stationary, have to decay in time. In section \ref{sec:consistentdiffusion}, a strategy is introduced that allows to construct schemes with a given finite difference approximation that is kept stationary. Such schemes must satisfy the additional requirement of stability under explicit time integration. No direct causal connection between the property of stationarity preservation and of (linear) stability of the scheme seems obvious, although both employ the common language of the Fourier transform. The study of possibly existing connections is subject of future work. 

An example of a linear stability analysis for schemes in two spatial dimensions can be found in \cite{morton01}.

\section{Numerical stationary states for the acoustic equations}
\label{sec:acoustic}

The condition of Theorem \ref{thm:detcondition} for stationarity preservation leads to a number of interesting additional results when applied to the acoustic equations. They appear as a linearization of the Euler equations and are used and studied in \cite{morton01}, \cite{amadori2015}, \cite{dellacherierieper10}, \cite{barsukow17} among others. They are of the form \eqref{eq:schemegeneral} and read
\begin{align}
 \del_t \vec v + \frac{1}{\epsilon^2}\nabla p &= 0 \label{eq:eqnacousticv}\\
 \del_t p + c^2 \, \div \vec v &= 0 \label{eq:eqnacousticp}
\end{align}
with $\vec v := (u, v)$, $q:= (u, v, p)$ in two spatial dimensions. An extension to three spatial dimensions is straightforward. 

Analogously to the low Mach number limit for the Euler equations (\cite{klainerman81}, \cite{klein95}, \cite{metivier01}) a $1/\epsilon^2$-scaling of the pressure gradient has been inserted (compare \cite{dellacherierieper10}). Formally, in the limit $\epsilon \to 0$ the solutions to \eqref{eq:eqnacousticv}--\eqref{eq:eqnacousticp} have constant pressure and a divergenceless velocity.

\begin{theorem} \label{thm:lowmachpdelevel}
The limit of low Mach number ${ \epsilon} \to 0$ for \eqref{eq:eqnacousticv}--\eqref{eq:eqnacousticp} is the same 
as the long time limit $t \to \infty$ of
\begin{align}
 \del_t \vec v + \nabla p &= 0\\
 \del_t p + c^2 \nabla \cdot \vec v &= 0
\end{align}
\end{theorem}
\begin{proof} Upon explicit calculation, the eigenvalues of $\vec J \cdot \vec k$ for the acoustic equations are $0$ and $\pm \frac{c |\vec k|}{\epsilon}$. Assembling the time evolution of the Fourier modes 
gives
 \begin{align}
  q(t, x) &= \int \dd \vec k \left( \hat q_0(\vec k) \exp(\ii \vec k \cdot \vec x) + \hat q_\pm(\vec k) \exp\left(\mp \ii |\vec k| \frac{c t}{\epsilon} + \ii \vec k \cdot \vec x\right) \right )
 \end{align}
 with $\hat q_0, \hat q_\pm$ following from the initial data. {Thus, if it exists, \begin{align}\lim_{\epsilon \to 0, \, t \text{ fixed}} \!\!q\,\,\,\, = \!\!\!\lim_{\epsilon\text{ fixed}, \,t \to  \infty} \!\!\!q\end{align}}
\end{proof}
Thus, decreasing $\epsilon$ by a factor of 10 and looking at the solution at time $t=1$, is the same as leaving $\epsilon$ as it was, and looking at the solution at time $t=10$. \new{The main result for the discrete situation, which establishes a connection between stationarity preservation and the low Mach number limit, is Theorem \ref{thm:lowmachdiscrete} below.}
 
Trivial stationary states of \eqref{eq:eqnacousticv}--\eqref{eq:eqnacousticp} are shear flows $\nabla p = 0$, $\del_x u = 0 =\del_y v$.

The acoustic equations are an example of systems considered in Theorem \ref{thm:stationarycontinuous}:
\begin{corollary} \label{cor:stationarypdelevelacoustic}
 The acoustic equations \eqref{eq:eqnacousticv}--\eqref{eq:eqnacousticp} allow for non-trivial stationary states given by
 \begin{align}
  \nabla p &= 0 & \nabla \cdot \vec v &= 0
 \end{align}
\end{corollary}
\begin{proof}
The eigenvector corresponding to the zero-eigenvalue of
\begin{align}
 \vec J \cdot \vec k &= \left( \begin{array}{ccc}  0&0&\frac{1}{\epsilon^2} k_x\\0&0&\frac{1}{\epsilon^2} k_y\\ c^2k_x & c^2	k_y & 0   \end{array} \right )
\end{align}
is
\begin{align}
 \hat q_0 = (-k_y, k_x, 0)^\text T \label{eq:acousticeigenvector}
\end{align}
to which $(\hat u, \hat v, \hat p)^\text T$ is only parallel if $\hat p = 0$ and $k_x \hat u + \hat v k_y = 0$, which are the Fourier transforms of $\nabla p = 0$, $\div \vec v = 0$. 
\end{proof}

From Theorem \ref{thm:continuousvorticity} follows the
\begin{corollary} \label{cor:vorticitypdelevelacoustic}
 The acoustic equations \eqref{eq:eqnacousticv}--\eqref{eq:eqnacousticp} have a constant of motion $\omega = \nabla \times \vec v$, i.e.
 \begin{align}
  \del_t \omega = 0 \label{eq:vorticityevolutionpdelevel}
 \end{align}
 Its Fourier transform is $-k_y \hat u + k_x \hat v$.
\end{corollary}
\begin{proof}
The left eigenvector corresponding to the zero eigenvalue of $\vec J \cdot \vec k$ is $(-k_y, k_x, 0)$ such that
\begin{align}
 (-k_y, k_x, 0) \del_t \veccc{\hat u}{\hat v}{\hat p} = 0
\end{align}
This amounts to a quantity $\omega$ whose Fourier transform is $k_x \hat v - k_y \hat u$, i.e. $\omega = \del_x v - \del_y u$, and $\del_t \omega = 0$.
\end{proof}
$\omega$ is called \emph{vorticity} by analogy with the Euler equations and it exemplifies Theorem \ref{thm:continuousvorticity}. Equation \eqref{eq:vorticityevolutionpdelevel} also immediately follows from the application of the curl operator to Equation \eqref{eq:eqnacousticv}. However the language of the Fourier transform is more useful when the discrete situation is considered.

\subsection{Vorticity preservation and the low Mach number limit}
\label{sec:vorticity}

The existence of a stationary vorticity for the acoustic equations has been given particular attention in \cite{torrilhon04}, \cite{mishra09preprint}, \cite{morton01} and others. The class of schemes that possess a stationary discrete counterpart has been given a name:

\begin{definition}
A consistent scheme for \eqref{eq:eqnacousticv}--\eqref{eq:eqnacousticp} is called \emph{vorticity-preserving}, if there exists a discretization of the vorticity that remains unchanged during the time evolution.
\end{definition}

By Theorem \ref{thm:discretevorticitygeneral} and Corollary \ref{cor:vorticitypdelevelacoustic}, having a vanishing eigenvalue of the evolution matrix yields the numerical stationary states as the right eigenvector and a numerical vorticity operator as the left eigenvector (as Fourier transforms). Conversely, if no eigenvalue of $\mathcal E(\vec k)$ vanishes for all $\vec k$, there is no discrete analogue of the vorticity that would remain stationary:

\begin{corollary} \label{eq:stationarityvorticityequiv}
 For the acoustic equations \eqref{eq:eqnacousticv}--\eqref{eq:eqnacousticp}, a scheme is vorticity preserving iff it is stationarity preserving.
\end{corollary}
\begin{proof} 
Theorem \ref{thm:discretevorticitygeneral}.
\end{proof}

It is known (see \cite{dellacherierieper10} among others) that the Roe scheme displays artefacts in the limit $\epsilon \to 0$ when applied to the equations \eqref{eq:eqnacousticv}--\eqref{eq:eqnacousticp}. With the concept of stationarity preservation these artefacts can be given a completely new interpretation.

Theorem \ref{thm:lowmachpdelevel} states that the limit $\epsilon \to 0$ of Equations \eqref{eq:eqnacousticv}--\eqref{eq:eqnacousticp} can be understood as the long time limit $t \mapsto \frac{t}{\epsilon} $, $\epsilon \to 0$ of
\begin{align}
 \del_t \vec v + \nabla p &= 0\\
 \del_t p + c^2 \nabla \cdot \vec v &= 0
\end{align}
The reason for that is the appearance of $\frac{c}{\epsilon}$ and $t$ only inside the combination $\frac{ct}{\epsilon}$. This is also true at discrete level:

\begin{theorem} \label{thm:lowmachdiscrete}
 Consider consistent and stable linear numerical schemes for \eqref{eq:eqnacousticv}--\eqref{eq:eqnacousticp}.
 Additionally, the eigenvalues of their evolution matrices $\mathcal E(\vec k)$ shall be linear in $c$.

 Vorticity preserving schemes that fulfill these conditions discretize the limit equations
 \begin{align}
  \nabla p &= 0 & \nabla \cdot \vec v &= 0
 \end{align}
 of \eqref{eq:eqnacousticv}--\eqref{eq:eqnacousticp} in the limit $\epsilon \to 0$.
\end{theorem}
\begin{proof}
 The proof consists of three parts:
 \begin{enumerate}[i)]
 \item There are asymptotic scalings (compare e.g. \cite{guillard99}, \cite{barsukow16}) of the dependent and independent variables $t$, $\vec x$, $p$, $\vec v$ and of $c$ which lead from
 \begin{align}
  \del_t \vec v + \nabla p &= 0\\
  \del_t \vec v + c^2 \nabla \cdot \vec v &= 0
 \end{align}
 to the rescaled equations
 \begin{align}
  \del_t \vec v + \frac{\nabla p}{\epsilon^2} &= 0\\
  \del_t \vec v + c^2 \nabla \cdot \vec v &= 0
 \end{align}
 These scalings can be easily computed explicitly and are, in their most general form,
 \begin{align}
  t &\mapsto \epsilon^\mathfrak a t & \vec v &\mapsto \epsilon^\mathfrak c \vec v \label{eq:rescaling1}\\
  \vec x &\mapsto \epsilon^\mathfrak b \vec x & p &\mapsto \epsilon^\mathfrak d p\\
  \vec k &\mapsto \epsilon^{-\mathfrak b} \vec k & c &\mapsto \epsilon^\mathfrak e c \label{eq:rescaling3}
 \end{align}
 where the free parameters $\mathfrak a, \mathfrak b, \mathfrak c, \mathfrak d, \mathfrak e$ have to satisfy
 \begin{align}
  \mathfrak d + \mathfrak a - \mathfrak b - \mathfrak c &= -2 \label{eq:rescaling4}\\
  2 \mathfrak e + \mathfrak c - \mathfrak b + \mathfrak a - \mathfrak d &= 0 \label{eq:rescaling5}
 \end{align}
 \item Now an analogous result to that of Theorem \ref{thm:lowmachpdelevel} is shown. The time evolution of the discrete Fourier modes, according to Equation \eqref{eq:semidiscrete} is given by the eigenvalues of the evolution matrix $\mathcal E(\vec k)$. 
 By assumption, they are linear in $c$.
 In general, the eigenvalue will have the form $c |\vec k| f(\vec k, \Delta x)$, with $f$ an arbitrary function (which does not depend on $c$).
 
 One now needs to show that the eigenvalues actually are linear in $\frac{c}{\epsilon}$. 
 
 An eigenvalue $c |\vec k| \cdot f$ of $\mathcal E(\vec k)$ means that there is a quantity $\hat q$ that satisfies an equation of the form 
 \begin{align}
 \del_t \hat q + c |\vec k| f \hat q &= 0 \label{eq:anyequationtoberescaled}
 \end{align}
 Rescaling according to \eqref{eq:rescaling1}--\eqref{eq:rescaling3}, and using \eqref{eq:rescaling4}--\eqref{eq:rescaling5} forces the rescaled form of Equation \eqref{eq:anyequationtoberescaled} to become 
 \begin{align}
  \del_t \hat q + \frac{c}{\epsilon} |\vec k| f \hat q &= 0
 \end{align}
 Therefore, the eigenvalues of $\mathcal E(\vec k)$ are linear in $\frac{c}{\epsilon}$. Therefore in the time evolution of any non-stationary Fourier mode only the combination $\frac{ct}{\epsilon}$ appears. This is in full analogy to the continuous case.
 
 \item In order to study the limit of low Mach numbers one looks at the long time evolution of the numerical scheme for the non-rescaled equations
 \begin{align}
  \del_t \vec v + \nabla p &= 0\\
  \del_t \vec v + c^2 \nabla \cdot \vec v &= 0
 \end{align}
 
 The statement of von Neumann stability is that every Fourier mode is either stationary, or decaying. Thus, as the numerical scheme is stable by assumption, then after long times only stationary Fourier modes will have survived. After long times, or equivalently for low Mach numbers, the numerical solution will be approaching one of the numerical stationary states of the scheme. As the scheme is assumed to be vorticity preserving, or equivalently stationarity preserving (by Corollary \ref{eq:stationarityvorticityequiv}), its stationary states are a discretization of the analytic stationary states (by Theorem \ref{thm:detcondition}), which by Theorem \ref{thm:lowmachpdelevel} are the limit equations for low Mach number.

 \end{enumerate}
\end{proof}

If the scheme is not stationarity preserving, it will fail to have the right limit as $\epsilon \to 0$, as its limit equations do not discretize all of the limit equations of the underlying PDE.

In fact, by considering the physical dimension that an eigenvalue of the evolution matrix $\mathcal E$ must have, one observes that linearity in $c$ is the only way for it to obtain the correct units. Therefore this assumption is actually always true.

This result is reminiscent of the concept of \emph{asymptotic preserving} (see e.g. \cite{jin99}). The general statement is in the following exemplified for a particular numerical setup (in Section \ref{sec:lowmach}).

\subsection{Notation}
\label{sec:notation}

In order to shorten the expressions, the following bracket notation is used:
\begin{align}
 [q]_{i+\frac12} &:= q_{i+1} - q_i & \{q\}_{i+\frac12} &:= q_{i+1} + q_i\\
 [q]_{i\pm1} &:= q_{i+1} - q_{i-1}\\
 [[q]]_{i\pm\frac12} &:= [q]_{i+\frac12} - [q]_{i-\frac12} & \{ \{ q \}\}_{i\pm\frac12}&:= \{q\}_{i+\frac12} + \{q\}_{i-\frac12}
\end{align}
The only nontrivial identity is
\begin{align}
 \{[q]\}_{i\pm\frac12} = [q]_{i+\frac12} + [q]_{i-\frac12} = [q]_{i\pm1}
\end{align}
For multiple dimensions the notation is combined, e.g.
\begin{align}
 [[q]_{i\pm1}]_{j\pm1} = q_{i+1,j+1} - q_{i-1,j+1} - q_{i+1, j-1} + q_{i-1,j-1}
\end{align}
The brackets for different directions commute.

\subsection{Dimensionally split schemes}
\label{sec:example-dimsplit}

For the purpose of the illustration of stationarity preserving properties of schemes in this section $\epsilon$ is treated as some finite parameter. The following section \ref{sec:lowmach} illustrates the limit $\epsilon \to 0$ in the context of stationarity preservation.

Consider a centered scheme with numerical diffusion for the acoustic system \eqref{eq:eqnacousticv}--\eqref{eq:eqnacousticp}, which has the general shape
\begin{align}
 \del_t q &+ \frac{1}{2 \Delta x} \Big ( J_x(q_{i+1,j} - q_{i-1, j}) - D_x(q_{i+1,j} - 2 q_{ij} + q_{i-1,j}) \Big ) \\
 &+ \frac{1}{2 \Delta y} \Big ( J_y(q_{i,j+1} - q_{i,j-1}) - D_y(q_{i,j+1} - 2 q_{ij} + q_{i,j-1}) \Big ) = 0 \label{eq:dimensionallysplitscheme}
\end{align}
A dimensionally split scheme typically has the following general form of the diffusion matrices $D_x$, $D_y$:
\begin{align}
 D_x &= \left ( \begin{array}{ccc} a_1 & 0 & a_2 \\ 0 & 0 & 0\\ a_3 & 0 & a_4 \end{array} \right ) &
 D_y &= \left ( \begin{array}{ccc} 0 & 0 & 0 \\ 0 & a_1 & a_2\\ 0 & a_3 & a_4 \end{array} \right ) \label{eq:diffusiondimsplit2d}
\end{align}

Dimensionally split schemes under certain conditions can be stationarity preserving:

\begin{theorem}[Stationarity preserving dimensionally split schemes] \label{thm:dimsplit}
 The dimensionally split scheme \eqref{eq:dimensionallysplitscheme} with \eqref{eq:diffusiondimsplit2d} is stationarity preserving if $a_1=0$. The stationary states fulfill $p = \const$ and
 \begin{align}
 \frac{[u]_{i\pm1,j} }{2 \Delta x} + \frac{[v]_{i,j\pm1} }{ 2\Delta y} - \frac{a_3}{c^2}\left(  \frac{ [[ u ]]_{i\pm\frac12,j}}{2\Delta x} + \frac{ [[v]]_{i,j\pm\frac12}}{2\Delta y} \right ) &= 0 \label{eq:discretedivcondition}
\end{align}
 which is a discretization of $\div \vec v = 0$.
\end{theorem}
\begin{proof}
The evolution matrix is easily found to be {\tiny
\begin{align}
 \mathcal E = \ii \left( \begin{array}{ccc} -\frac{ a_1 (t_x - 2 + \frac{1}{t_x})}{2\Delta x} & 0 & - \frac{ a_2 (t_x - 2 + \frac{1}{t_x})}{2\Delta x} + \frac{ (t_x - \frac{1}{t_x})}{2 \Delta x \epsilon^2}\\
                      0 & -\frac{ a_1 (t_y - 2 + \frac{1}{t_y})}{2\Delta y} & - \frac{ a_2 (t_y - 2 + \frac{1}{t_y})}{2\Delta y} + \frac{ (t_y - \frac{1}{t_y})}{2\Delta y \epsilon^2}\\
		- \frac{ a_3 (t_x - 2 + \frac{1}{t_x})}{2\Delta x} + \frac{c^2 (t_x - \frac{1}{t_x}) }{2 \Delta x} & - \frac{ a_3 (t_y - 2 + \frac{1}{t_y})}{2\Delta y} + \frac{c^2  (t_y - \frac{1}{t_y}) }{2 \Delta y} & - \frac{ a_4 (t_x - 2 + \frac{1}{t_x})}{ 2\Delta x} - \frac{ a_4 (t_y - 2 + \frac{1}{t_y}) }{2 \Delta y}
                     \end{array} \right ) \label{eq:evolutionmatrix}
\end{align}}
whose determinant is only zero (independently of $\vec k$), if $a_1 = 0$ as can be shown upon direct computation. In this case the corresponding eigenvector is
\begin{align}
 \veccc{\frac{a_3 (t_y - 2 + t_y^{-1})}{2\Delta y} - \frac{c^2  (t_y - t_y^{-1}) }{2 \Delta y}}{- \frac{ a_3 (t_x - 2 + t_x^{-1})}{2\Delta x} + \frac{c^2  (t_x - t_x^{-1}) }{2 \Delta x} }{0} \label{eq:dimspliteigenvector}
\end{align}
which amounts, by inverting the Fourier transform, to the given discrete divergence operator. 
\end{proof}

Numerical data that exactly satisfy \eqref{eq:discretedivcondition} remain unchanged during the evolution. The discrete operator \eqref{eq:discretedivcondition} is a first order discretization of $\del_x u + \del_y v$, if $a_3 \neq 0$. Choosing both $a_1 = 0$ and $a_3 = 0$ in \eqref{eq:discretedivcondition} and \eqref{eq:diffusiondimsplit2d} makes all the spatial operators reduce to central differences:

\begin{corollary} \label{thm:central}
 A scheme for the system \eqref{eq:eqnacousticv}--\eqref{eq:eqnacousticp} whose spatial derivatives are discretized by central differences in two spatial dimensions is stationarity preserving.
\end{corollary}

Choosing $a_3 = 0$ means that the discrete divergence operator which is exactly preserved during the time evolution, is a central one. Together with $a_1 = 0$, however, this would mean that there is no diffusion on the velocity variables at all. In practice this is often not desirable as then the scheme will not be stable upon usage of an explicit time integrator (e.g. forward Euler). One might wonder whether there exists a discrete velocity diffusion such that the resulting scheme would keep the central divergence exactly stationary. This would lead to the stationary states being discretized to higher order. One is thus led to the more general question of finding a discrete diffusion that vanishes whenever a given discrete divergence does. This is discussed in section \ref{sec:consistentdiffusion}.

The upwind, or Roe, scheme has $D_x = |J_x|$, $D_y = |J_y|$, with the absolute value being defined on the eigenvalues. This gives
\begin{align}
 D_x &= \left(  \begin{array}{ccc} \frac{c}{\epsilon} \\&0\\&& \frac{c}{\epsilon} \end{array} \right ) &
 D_y &= \left(  \begin{array}{ccc} 0 \\&\frac{c}{\epsilon}\\&& \frac{c}{\epsilon} \end{array} \right )
\end{align}
which is of the form \eqref{eq:diffusiondimsplit2d}, but violates the condition $a_1 = 0$ found in Theorem \ref{thm:dimsplit}. Thus one has proved the

\begin{corollary}
 The Roe scheme for the system \eqref{eq:eqnacousticv}--\eqref{eq:eqnacousticp} in two spatial dimensions is not stationarity preserving.
\end{corollary}

\begin{figure}[h]
 \centering
 \includegraphics[width=0.48\textwidth]{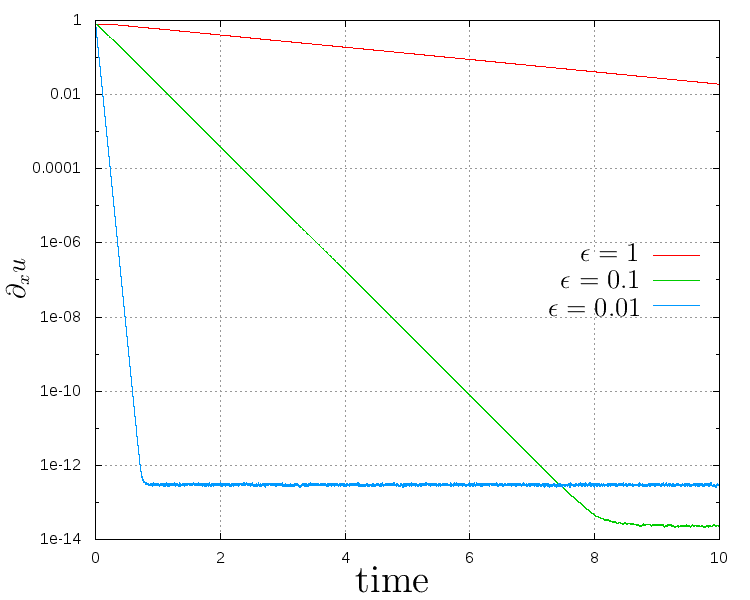}
 \hfill \includegraphics[width=0.48\textwidth]{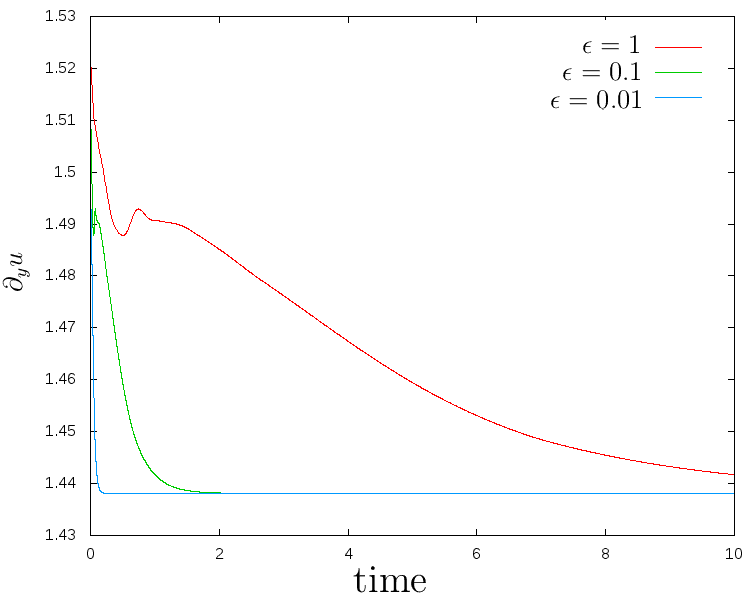}
 \caption{{\sl Left}: Decay of $\displaystyle\|\del_x u\|_{L^1} \sim \exp\left( -\frac{t}{\epsilon}  \right )$. {\sl Right}: $\|\del_y u\|_{L^1}$. Both figures were measured for a stationary vortex setup in a simulation using the Roe solver for the acoustic equations and show curves for values $\epsilon = 1, 0.1, 0.01$. Note the very different vertical axis scalings in the two plots: whereas $\del_x u \neq 0$ does not comply with stationarity for the Roe scheme, after some transients the scheme settles down on a shear flow ($\del_y u \neq 0$) that is not significantly different from the initial data.}
 \label{fig:greshodecay}
\end{figure}

\begin{figure}[h]
 \centering
 \includegraphics[width=0.3\textwidth]{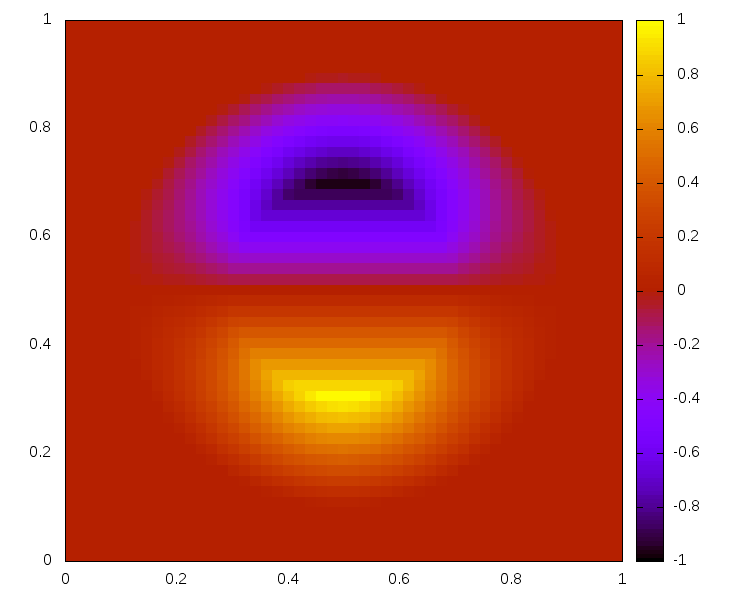}
 \includegraphics[width=0.3\textwidth]{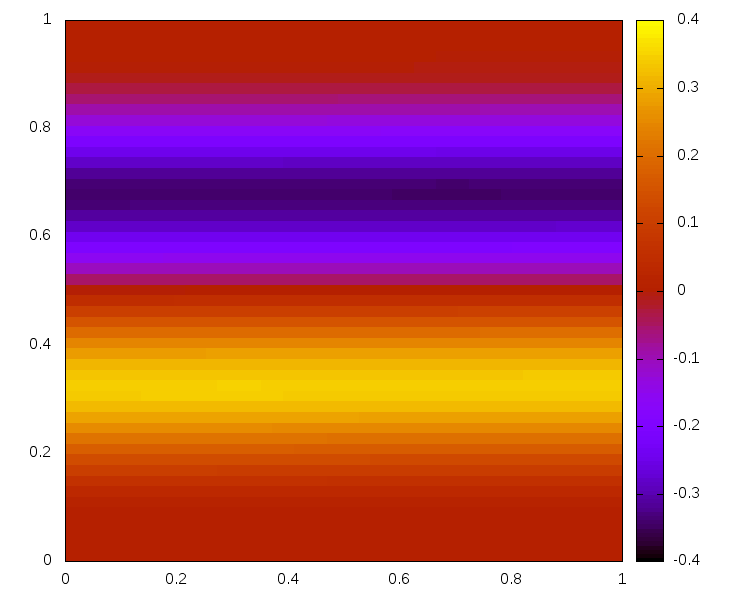}\\
 \includegraphics[width=0.3\textwidth]{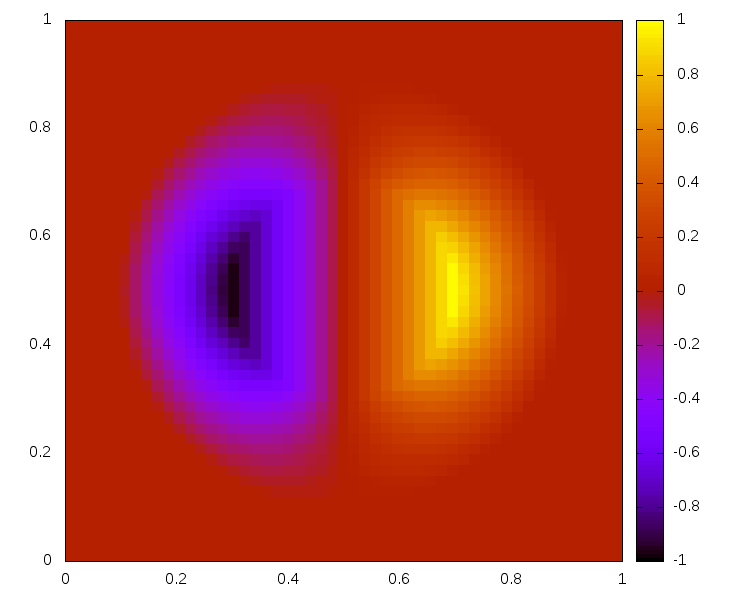}
 \includegraphics[width=0.3\textwidth]{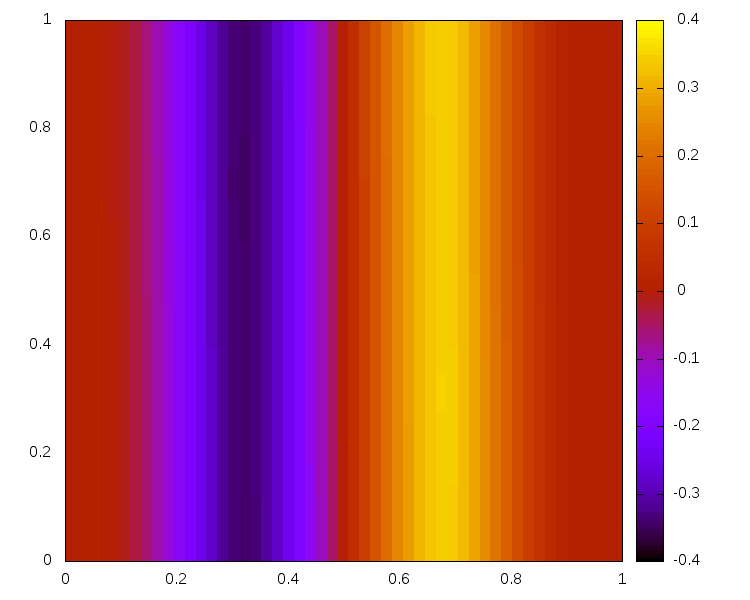}
 \caption{Simulation results for $\epsilon = 10^{-3}$ of a vortex setup with the Roe scheme. {\sl Left} are results at time $t=0$, {\sl right} -- at $t=0.3$. {\sl Top row}: $u$ (color coded), {\sl bottom row}: $v$ (color coded). The Roe scheme fails to keep the setup stationary and transitions to a trivial stationary state (shear flow).}
 \label{fig:greshoresults}
\end{figure}

This can be observed in the experiment. As an example the numerical time evolution of a stationary divergence-free vortex flow around the origin is shown in Figures \ref{fig:greshodecay}--\ref{fig:greshoresults}. From an initial state (that is derived from the analytic stationary solution) one observes the numerical solution to move over to some other stationary solution. Initially, $\del_x u + \del_y v = 0$, but $\del_x u \neq 0$ in general. One observes however in Figure \ref{fig:greshodecay} that the Roe scheme diffuses away $\del_x u$ exponentially in time, until it reaches values comparable with machine precision. The initial velocities, shown in the left column of Fig. \ref{fig:greshoresults}, are modified such that after long times only a shear flow is left over. The only states that the scheme is able to keep stationary, are trivial ones. By stability, the scheme is diffusing away all the others. On the other hand, a stationarity preserving scheme would keep stationary also discrete versions of vortical, and in general of all divergenceless flows. 

The behaviour of the scheme as $\epsilon \to 0$ is discussed in section \ref{sec:lowmach}.
 
For the above example the discrete vorticity that is preserved exactly during the time evolution is, by \eqref{eq:dimspliteigenvector} (the notation having been introduced in section \ref{sec:notation})
\begin{align}
 \frac{[v]_{i\pm1,j} }{2 \Delta x} - \frac{[u]_{i,j\pm1} }{ 2\Delta y} + \frac{a_3}{c^2}\left( \frac{ [[u]]_{i,j\pm\frac12}}{2\Delta y} - \frac{ [[ v ]]_{i\pm\frac12,j}}{2\Delta x} \right ) = \del_x v - \del_y u + \mathcal O(\Delta x, \Delta y)
\end{align}

In \cite{morton01} the Appendix deals with the preservation of a specific discrete vorticity for a certain family of schemes, and with the production rate in case of non-preservation. However this analysis, as well as the treatment of the acoustic equations in \cite{lerat07}, \emph{assume} a vorticity discretization. There might still exist some other vorticity discretization that is exactly preserved. Therefore a more adequate procedure would be to first check (via the determinant of the evolution matrix) the existence of \emph{any} preserved discrete vorticity and to find its shape by evaluating the eigenvector corresponding to the vanishing eigenvalue. Only if the evolution matrix does not contain any vanishing eigenvalues can one claim that there is no preserved discrete vorticity. The condition of Theorem 2 formulated in \cite{lerat07} therefore is sufficient, but not necessary for vorticity preservation.

\subsection{Low Mach number limit}
\label{sec:lowmach}

Next the vortex setup as in Figures \ref{fig:greshodecay}--\ref{fig:greshoresults} is considered for the limit of $\epsilon \to 0$. It is solved with the upwind/Roe scheme, which has only trivial stationary states.

The initial data of the vortex are a discrete version of an analytically stationary solution and have been thus obtained from a divergence-free solution. The upwind/Roe scheme, since it is stable, keeps certain states exactly stationary and diffuses everything else away with time. This diffusion time scales with $\epsilon$, because the non-zero eigenvalues of the evolution matrix scale with $1/\epsilon$. After long time therefore one is left with a numerical stationary state of the scheme. 

If the set of numerical stationary states, however, consists only of trivial ones
(as it is the case for the Roe scheme), then this numerical solution will have lost all resemblance to the analytic one. In this example the vortex is diffused away and a shear flow is left over. Therefore the observed ``low Mach number artefacts'' are entirely due to the scheme's stationary solutions not being discretizations of all the analytic ones.

The low Mach number limit $\epsilon \to 0$ for the acoustic equations makes the scheme attain a numerical stationary state on time scales $\mathcal O(\epsilon)$. In order to improve the quality of the numerical solution therefore one needs to choose a scheme which has nontrivial numerical stationary states that capture the rich set of nontrivial stationary states, i.e. a stationarity preserving scheme.

Reducing $\Delta x$ does not really solve the problem, because then the diffusion time will be longer, but the stationary states will still not be any more similar to the analytic ones. In the literature, there already exist several strategies that have been developed in order to cope with the low Mach number problems. They can now be understood in the light of the new arguments that employ the idea of a stationarity preserving scheme. Adapting the matrices to the case of acoustic equations yields the following selection of diffusion matrices:

\begin{enumerate}[1.]
 \item Method from \cite{barsukow16}: $\displaystyle D_x = \left ( \begin{array}{ccc} 0 &0 & \frac{1}{\epsilon^2} \\0&0&0\\ -c^2 &0& 0 \end{array} \right ) $
 \item Method from \cite{dellacherierieper10}: $\displaystyle D_x = \left ( \begin{array}{ccc} 0 &0& 0 \\0&0&0\\ -c^2 &0& \frac{2c}{\epsilon} \end{array} \right ) $
 \item A new method: $\displaystyle D_x = \left ( \begin{array}{ccc} 0 &0& \frac{1}{\epsilon^2} \\0&0&0\\ 0 &0& \frac{2c}{\epsilon} \end{array} \right ) $ \label{it:newsolver}
\end{enumerate}

\begin{figure}[h]
 \centering
 \includegraphics[width=0.68\textwidth]{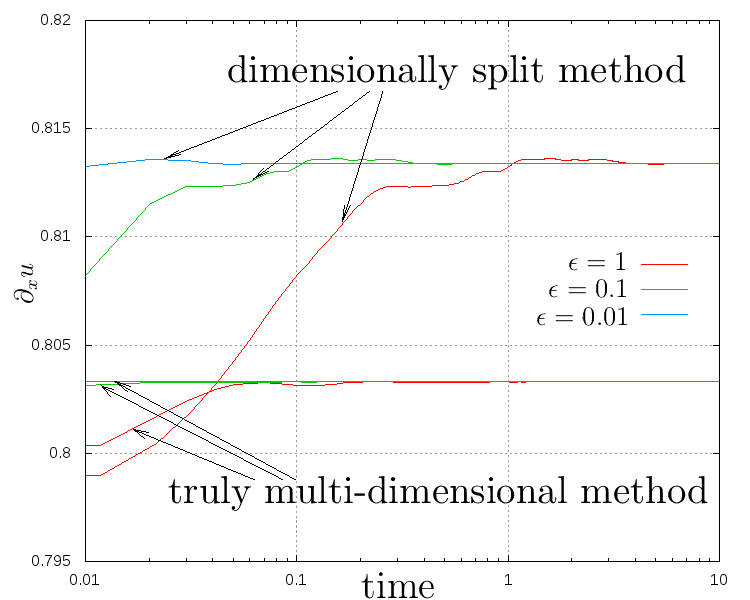}
 \caption{Time evolution of $\|\displaystyle\del_x u\|_{L^1}$. Figure measured for a stationary vortex setup in a simulation using solver \ref{it:newsolver} mentioned above (``dimensionally split method'') and the multidimensional solver (``truly multi-dimensional method'') presented in section \ref{sec:consistentdiffusion} for the acoustic equations and show curves for values $\epsilon = 1, 0.1, 0.01$. (Note the scalings of both axes.) Observe the absence of diffusion (contrary to Fig. \ref{fig:greshodecay}, left) and the improved quality of the simulation upon usage of a truly multi-dimensional solver.}
 \label{fig:greshostationary}
\end{figure}

\begin{figure}[h]
 \centering
 \includegraphics[width=0.48\textwidth]{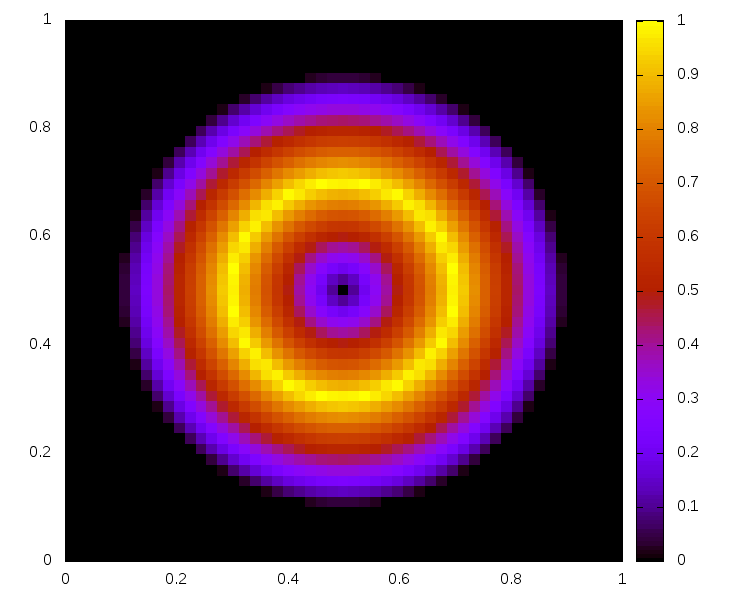}
 \includegraphics[width=0.48\textwidth]{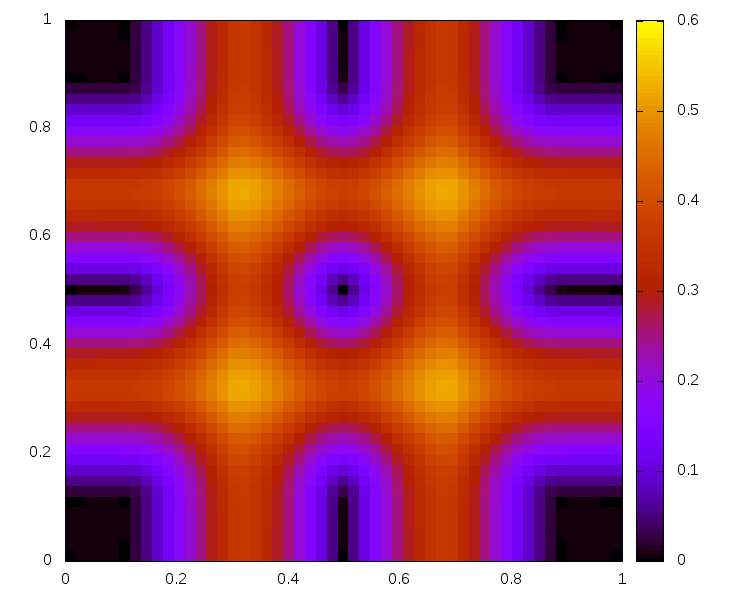}\\
 \includegraphics[width=0.48\textwidth]{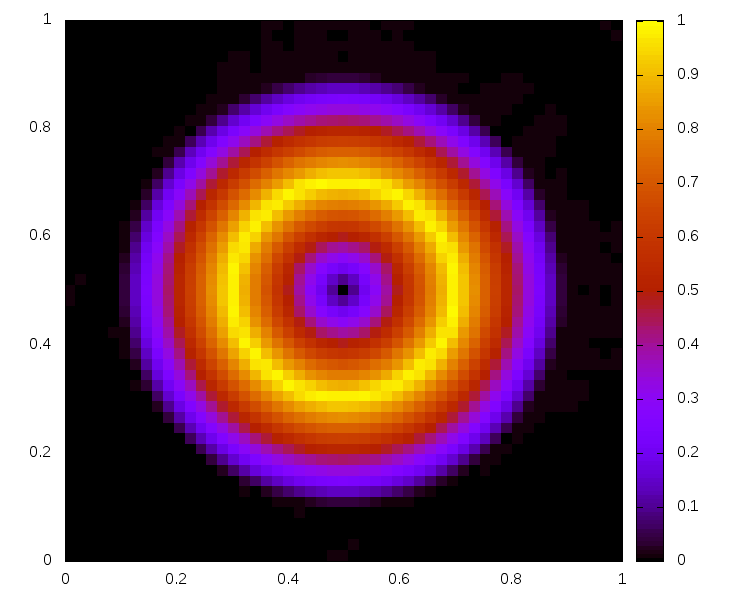}
 \hfill \includegraphics[width=0.48\textwidth]{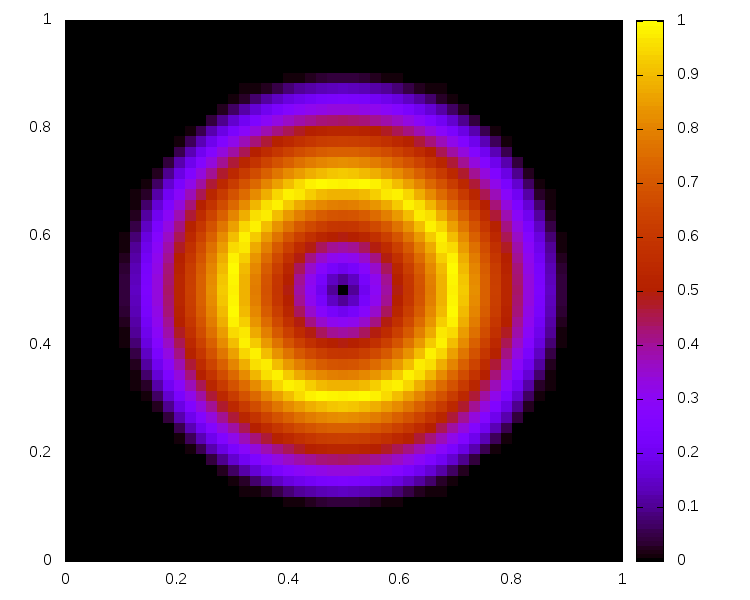}
 \caption{{\sl Top left}: Initial setup of a stationary divergence-free vortex. {\sl Top right}: Solution at $t=10$ with the upwind (Roe) scheme. {\sl Bottom left}: Solution at $t=10$ with solver \ref{it:newsolver} from the above list. {\sl Bottom right}: Solution at $t=10$ with the truly multi-dimensional solver discussed in section \ref{sec:consistentdiffusion}. $\sqrt{u^2 + v^2}$ is color coded; all simulations performed on a $50 \times 50$ grid with forward Euler; all methods are of first order in space and time. Observe the improved quality when using a stationarity preserving scheme.}
 \label{fig:greshostationaryfinal}
\end{figure}

Note how all of them have $a_1 = 0$, and are thus stationarity preserving. \new{Experimental results obtained with method no. \ref{it:newsolver} of the above list are shown in Fig. \ref{fig:greshostationary}--\ref{fig:greshostationaryfinal}. This method has been found to be stable under explicit time integration experimentally.}


Methods that use implicit time integration solve an additional problem: in the limit $\epsilon \to 0$ the time step of an explicit time integrator that solves \eqref{eq:eqnacousticv}--\eqref{eq:eqnacousticp} must scale as $\mathcal O(\epsilon)$ if it shall fulfill the CFL condition. This restriction is not present for an implicit solver. However, in order to be able to resolve the limit $\epsilon \to 0$ numerically, the discretization of the spatial derivatives has to fulfill the same conditions, whether the time integration is explicit or implicit. Such solvers (e.g. \cite{cordier12}, \cite{reiss14} and many others) often use central differences, and thus choose a very easy way to obtain a stationarity preserving scheme (as central differences are stationarity preserving by the Corollary \ref{thm:central} from section \ref{sec:example-dimsplit}).

\section{Constructing stationarity preserving schemes: consistent diffusion}
\label{sec:consistentdiffusion}

Central differences for the spatial derivatives have been shown to be stationarity preserving in Corollary \ref{thm:central}, but they are known to be unstable under forward Euler time integration. Therefore the natural question is whether it is possible to write down a stabilizing diffusion for them that does not spoil the property of stationarity preservation. This exemplifies how the theory developed so far can be used constructively.

\subsection{Continuous case}

Consider, as an example, the linear system \eqref{eq:schemegeneral} in $d=2$ spatial dimensions
\begin{align}
 \del_t q + J_x \del_x q + J_y \del_y q &= 0 \label{eq:schemegeneral2d}
\end{align}
with $q : \mathbb R^+_0 \times \mathbb R^2 \to \mathbb R^n$ and $J_x$, $J_y$ being $n \times n$ matrices. The stationary solutions are given by
\begin{align}
  J_x \del_x q + J_y \del_y q &= 0 \label{eq:schemegeneralstationary2d}
\end{align}
Consider now a numerical scheme for Equation \eqref{eq:schemegeneral2d}, e.g. a finite volume scheme or a finite difference scheme. Before the concepts can be detailed for the discrete case, it is easier to discuss them in a continuous situation (as is done e.g. in \cite{sidilkover02}, \cite{lerat07}), where effects of the numerics are taken into account as a diffusive term, e.g. as
\begin{align}
 \del_t q + J_x \del_x q + J_y \del_y q &= D_x \del_x^2 q + D_y \del_y^2 q + \bar D \del_x \del_y q \label{eq:schemegeneraldiffusion2d}
\end{align}
with $D_x$, $D_y$, $\bar D$ matrices. Of course, for stability there are certain conditions that these matrices need to fulfill, which shall not matter for the moment.

Consider now initial data fulfilling \eqref{eq:schemegeneralstationary2d}, such that they are preserved exactly in time if the evolution is governed by \eqref{eq:schemegeneral2d}. If the initial data are evolved according to \eqref{eq:schemegeneraldiffusion2d}, then their initial evolution $\del_t q$ is given entirely by the diffusion
\begin{align}
 D_x \del_x^2 q + D_y \del_y^2 q + \bar D \del_x \del_y q
\end{align}
In this situation, a {\sl stationarity consistent diffusion} would be a term containing second derivatives, that vanishes whenever $J_x \del_x q + J_y \del_y q$ vanishes.

As an example, one could take 
\begin{align}
J_x \del_x^2 q + J_y \del_x \del_y q = \del_x (J_x \del_x q + J_y \del_y q) \label{eq:consistentdiffusioncontinuous}
\end{align}
i.e.
\begin{align}
 D_x &= J_x & D_y &= 0 & \bar D = J_y
\end{align}
or anything proportional to it. Observe also the necessary appearance of mixed second derivatives.

To implement this idea in a fully discrete situation is subject of the rest of this paper.

\subsection{Discrete case}

\begin{definition}
 A finite difference formula $\mathcal B$, that vanishes whenever another (given) finite difference $\mathcal A$ does, shall be called \emph{stationarity-consistent} with $\mathcal A$.
\end{definition}

The vanishing of a finite difference formula is studied via the Fourier transform again. $\mathcal B$ being stationarity-consistent with $\mathcal A$ means that the Fourier transform of $\mathcal B$ contains the Fourier transform of $\mathcal A$ as a factor. Clearly, if $\mathcal A$ is a stationarity preserving discretization of $\vec J \cdot \nabla q$, then adding a stationarity-consistent diffusion does not destroy the stationarity preservation property.

Note that according to \eqref{eq:fouriertrafostencilcomplete} the Fourier transform of any compact linear finite difference $\sum_S \alpha_S q_{I + S}$ is proportional to the Laurent polynomial
\begin{align}
 \sum_S \alpha_S \hat q \prod_{m=1}^d t_m^{S_m}
\end{align}
in the variables $\{ t_m \}_{1 \leq m \leq d}$. Note that $\hat q$ appears linearly. This establishes a mapping between compact linear finite difference formulae in $d$ dimensions and Laurent polynomials in $d$ variables. The condition of stationarity-consistency thus means that the polynomial corresponding to the linear finite difference $\mathcal B$ can be written as the polynomial of the linear finite difference $\mathcal A$ times some factor; this factor however does not need to by a polynomial itself. Note that again, in both $\mathcal A$ and $\mathcal B$, $\hat q$ is assumed to appear linearly. Therefore, for scalar problems, any two finite difference formulae are stationarity-consistent. This is different when systems are considered, as is exemplified in the following.

\subsection{Diffusion consistent with a stationary divergence}

For the rest of this section the focus is again on the acoustic equations \eqref{eq:eqnacousticv}--\eqref{eq:eqnacousticp} in two spatial dimensions, and in particular on the divergence operator appearing therein. Having first order schemes in mind, the stencils involve only the cell itself and its eight neighbours.

\begin{definition}
 The \emph{Moore neighborhood} of cell $(i,j)$ is the cell itself and the 8 cells $\{(i, j \pm 1)\} \cup \{ (i \pm 1, j) \} \cup \{ (i \pm 1, j \pm 1) \}$.\\
 A \emph{Moore stencil} at cell $(i,j)$ is a stencil involving only cells from the Moore neighborhood of cell $(i,j)$.
\end{definition}

Consider a central divergence approximation in two spatial dimensions
\begin{align}
 \frac{[u]_{i\pm1,j} }{2 \Delta x} + \frac{[v]_{i,j\pm1} }{ 2\Delta y} \label{eq:centraldiv}
\end{align}

\begin{theorem} \label{thm:nocentraldiv}
 There is no non-zero compact linear finite difference formula on a Moore stencil that approximates second derivatives of $u$ and $v$ and which is stationarity-consistent with the central divergence.
\end{theorem}

\textit{Note}. This explains why the authors in \cite{jeltsch06}, \cite{torrilhon04} found that ``the choice of central differences turned out to be not very fruitful''.

\begin{proof}
The Fourier transform of \eqref{eq:centraldiv} is, up to an irrelevant prefactor,
\begin{align}
 \hat u \frac{t_x - t_x^{-1} }{2 \Delta x} + \hat v \frac{t_y - t_y^{-1} }{ 2\Delta y} \label{eq:centraldivfourier}
\end{align}
This is a Laurent polynomial in $t_x$ and $t_y$. On a Moore neighborhood, the maximal and minimal powers of $t_x$ that can appear, are $t_x$ and $\frac{1}{t_x}$, analogously for $t_y$.

In order to find a stationarity consistent diffusion one needs to find a factor that would make \eqref{eq:centraldivfourier} a discretization of second derivatives, e.g. a discretization of \eqref{eq:consistentdiffusioncontinuous}. The Fourier transform of such a finite difference formula must contain a prefactor of $k_x^2$ and thus of $(t_x-1)^2$. Multiplying \eqref{eq:centraldivfourier} with $(t_x-1)$ yields
\begin{align}
 \hat u \frac{(t_x - 1)^2(t_x + 1) }{2 \Delta x t_x} + \hat v \frac{(t_y - 1)(t_x - 1)(t_y + 1) }{ 2\Delta y t_y} \label{eq:halfsecondderivfourier}
\end{align}
However now the first term cannot be constructed on a Moore stencil. Therefore the correct factor would have rather been $\frac{t_x-1}{t_x + 1}$. (To divide out $(t_x+1)$ is the only option because the result has to be a Laurent polynomial again, i.e. one cannot divide by, say, $t_x+2$). However the division by $t_x+1$, possible in the first term, is impossible in the second term. 
\end{proof}

Modifying the divergence discretization allows to find a stationarity-consistent diffusion. Obviously, in order to be able to divide by $1+t_x$, the discrete divergence must contain this factor in all terms from the beginning. Analogously one needs to be able to divide by $1 + t_y$.

\begin{theorem}
The only symmetric linear divergence discretization on $3 \times 3$ cells that allows for a non-zero stationarity-consistent diffusion is
\begin{align}
  \frac{\{\{[u]_{i\pm1}\}\}_{j\pm\frac12} }{8 \Delta x } + \frac{[\{\{v\}\}_{i\pm\frac12}]_{j\pm1} }{ 8\Delta y} \label{eq:averageddivergence}
\end{align}
A linear stationarity-consistent diffusion discretization associated to the divergence \eqref{eq:averageddivergence} is  
\begin{align}
  &\frac14 c_1 \left( \frac{\{\{[[u]]_{i\pm\frac12} \}\}_{j\pm\frac12}}{\Delta x}+  \frac{ [[v]_{i\pm1}]_{j\pm1} }{\Delta y} \right )
    + \frac14 c_2 \left( \frac{     [[u]_{i\pm1}]_{j\pm1} }{\Delta x} + \frac{ [[\{\{v\}\}_{i\pm\frac12}]]_{j\pm\frac12}  }{\Delta y} \right ) \label{eq:secondderivfouriermultid}
\end{align}
with arbitrary parameters $c_1$, $c_2$.
\end{theorem}
\begin{proof}
{This constructive proof builds on the argumentation in the proof of Theorem \ref{thm:nocentraldiv}. There, the crucial observation was that the term
\begin{align}
 \hat v \frac{(t_y - 1)(t_x - 1)(t_y + 1) }{ 2\Delta y t_y} 
\end{align}
in \eqref{eq:halfsecondderivfourier} does not contain a factor $(t_x+1)$, which makes it impossible to construct a stationarity consistent diffusion. The Fourier transform of the discrete divergence can be made to contain such a term from the beginning. In order to make the corresponding finite difference formula symmetric, along with $(t_x + 1)$ one also has to include $(1 + \frac{1}{t_x})$. Thus \eqref{eq:centraldivfourier} is to be replaced by 
\begin{align}
 \hat u \frac{(t_x - 1)(t_x + 1)(t_y + 1)^2 }{8 \Delta x t_x t_y} + \hat v \frac{(t_y - 1)(t_y + 1)(t_x + 1)^2 }{ 8\Delta y t_y t_x} \label{eq:divfouriercompletemultid}
\end{align}
This is the Fourier transform of \eqref{eq:averageddivergence}.
The strategy presented in the proof of Theorem \ref{thm:nocentraldiv} can now be successfully implemented: \eqref{eq:divfouriercompletemultid} can now be multiplied with $\frac{t_x - 1}{t_x + 1}$ to yield the Laurent polynomial
\begin{align}
 \hat u \frac{(t_x - 1)^2(t_y + 1)^2 }{8 \Delta x t_x t_y} + \hat v \frac{(t_y - 1)(t_y + 1)(t_x - 1)(t_x + 1) }{ 8\Delta y t_y t_x} 
\end{align}
Taking the inverse Fourier transform\footnote{This is easy, as any term $\hat u \, t_x^a t_y^b$ has just to be replaced by $u_{i+a, j+b}$, and analogously for $v$.} yields the first term in \eqref{eq:secondderivfouriermultid}, and analogously the others can be obtained. \eqref{eq:secondderivfouriermultid} is a discrete second derivative because to highest order in an expansion in powers of $\Delta x$ it equals to
\begin{align}
 \frac12 c_1 &\Delta x (\del_x^2 u + \del_x\del_y v) + \frac12 c_2 \Delta y (\del_x \del_y u + \del_y^2 v)
 + \mathcal O(\Delta x^2, \Delta y^2)
\end{align}}
\end{proof}

The divergence discretization \eqref{eq:averageddivergence} is -- by equivalence of stationarity and vorticity preservation -- the ``extended operator'' in \cite{torrilhon04}, \cite{jeltsch06} and has also has been suggested in \cite{mishra09preprint} for the system wave equation. The above proof shows that there is actually no other choice among directionally unbiased finite difference formulae defined on a $3 \times 3$ grid, i.e. among symmetric Moore discretizations. In \cite{sidilkover02} some non-standard finite difference methods are introduced in the context of steady Euler equations. They are reminiscent of the stencils above for the linearized Euler equations.

Note that still it is possible that the schemes that contain this finite difference formula differ -- they might treat the pressure differently, or have different order (e.g. compare the scheme in \cite{morton01} to \cite{jeltsch06}). In order to complete the scheme one can use the same discrete operators for the spatial derivatives of the pressure and include a diffusion on the pressure that mimics the Laplacian $\div \,\grad p = \del_x^2 p + \del_y^2 p$. This gives a method that is very similar to the one suggested in \cite{mishra09preprint} and \cite{jeltsch06}:

{\tiny
\begin{align*}
   \del_t \veccc{u}{v}{p} &+ \frac{1}{8\Delta x}\left( \veccc{\frac{1}{{ \epsilon^2}} [ \{\{ p \}\}_{ j\pm\frac12}  ]_{ i\pm 1} }{0}{c^2 [ \{\{ u \}\}_{ j\pm\frac12}  ]_{ i\pm 1}} - \veccc{\frac{c}{ \epsilon} [[ \{\{ u \}\}_{ j\pm\frac12}]]_{ i\pm\frac12}  }{0}{\frac{c}{ \epsilon} [[ \{\{ p \}\}_{ j\pm\frac12}]]_{ i\pm\frac12}  } -  \veccc{0\phantom{]_{i+\frac12}}}{\frac{c}{ \epsilon} [ [u^{ n}]_{ j\pm1,.}]_{ i\pm1}  }{0\phantom{]_{i+\frac12}}}    \right )\\
   &+ \frac{1}{8\Delta y} \left( \veccc{0}{\frac{1}{{ \epsilon^2}} [ \{\{ p \}\}_{ i\pm\frac12}  ]_{ j\pm 1} }{c^2 [ \{\{ v \}\}_{ i\pm\frac12}  ]_{ j\pm 1}} - \veccc{0}{\frac{c}{ \epsilon} [[ \{\{ v \}\}_{ i\pm\frac12}]]_{ j\pm\frac12}  }{\frac{c}{ \epsilon} [[ \{\{ p \}\}_{ i\pm\frac12}]]_{ j\pm\frac12}  }  -  \veccc{\frac{c}{ \epsilon} [ [v^{ n}]_{ i\pm1,.}]_{ j\pm1}  }{0\phantom{]_{i+\frac12}}}{0\phantom{]_{i+\frac12}}}   \right ) = 0
\end{align*}}
This scheme reduces to the upwind (Roe) scheme if restricted to one spatial dimension. Experimentally, it shows stability up to a CFL number of 1 (rather than 0.5 that is found for dimensionally split schemes in two dimensions). It is stationarity preserving by construction. Results of a simulation of a divergence-free vortex can be seen in Fig. \ref{fig:greshostationary}--\ref{fig:greshostationaryfinal}, and there is evidence for a slight superiority of results obtained with this multi-dimensional scheme as compared to the dimensionally split method presented in the same figures.

Consider the first term in \eqref{eq:secondderivfouriermultid} in more detail (taking $\Delta y = \Delta x$ to ease the notation):

{\scriptsize 
\begin{align}
  &\frac14 \left( \frac{\{\{[[u]]_{i\pm\frac12} \}\}_{j\pm\frac12}}{\Delta x}+  \frac{ [[v]_{i\pm1}]_{j\pm1} }{\Delta y} \right ) 
  =  (\del_x^2 u +\del_x \del_y v) \Delta x
  \\&+\frac{ (3 \del_x^2 \del_y^2 u+\del_x^4 u+2 (\del_x \del_y^3 v+\del_x^3 \del_y v)) \Delta x^3}{12 }\\
  &+\frac{ (15 \del_x^2 \del_y^4 u+15 \del_x^4 \del_y^2 u+2 \del_x^6 u+6 \del_x \del_y^5 v+20 \del_x^3 \del_y^3 v+6 \del_x^5 \del_y v) \Delta x^5}{720 }\\
  &+\frac{ (14 \del_x^2 \del_y^6 u+35 \del_x^4 \del_y^4 u+14 \del_x^6 \del_y^2 u+\del_x^8 u+4 (\del_x \del_y^7 v+7 (\del_x^3 \del_y^5 v+\del_x^5 \del_y^3 v)+\del_x^7 \del_y v)) \Delta x^7}{20160}\\
  &+\mathcal O(\Delta x^9)
\end{align}}
Obviously, the highest order term is just of the form \eqref{eq:consistentdiffusioncontinuous} mentioned in the Introduction. It would however be erroneous to expect the higher order terms to be higher derivatives of $\del_x u + \del_y v$, because it is not $\del_x u + \del_y v$ that is preserved exactly by the numerics -- it is rather the discrete divergence
{\scriptsize
\begin{align}
   \mathscr D &= \frac{\{\{[u]_{i\pm1}\}\}_{j\pm\frac12} }{8 \Delta x } + \frac{[\{\{v\}\}_{i\pm\frac12}]_{j\pm1} }{ 8\Delta y} 
   = (\del_x u+\del_y v)
   \\&+\frac{1}{12} \Big (3 \del_x\del_y^2 u+2 (\del_x^3 u+ \del_y^3 v)+3 \del_x^2 \del_y v \Big ) \Delta x^2\\
   &+\frac{1}{240} \Big (5 \del_x\del_y^4 u+2 (5 \del_x^3 \del_y^2 u+\del_x^5 u + \del_y^5 v + 5 \del_x^2 \del_y^3 v)+5 \del_x^4 \del_y v \Big ) \Delta x^4\\
   &+\frac{  1  }{{10080}} \Big (7 \del_x \del_y^6 u+35 \del_x^3 \del_y^4 u+21 \del_x^5 \del_y^2 u+2 \del_x^7 u+2 \del_y^7 v +7 (3 \del_x^2 \del_y^5 v+5 \del_x^4 \del_y^3 v+\del_x^6 \del_y v) \Big ) \Delta x^6 \\ 
   &+\mathcal O(\Delta x^8)
\end{align}}

Indeed, one then finds that
\begin{align}
  &\frac14 \left( \frac{\{\{[[u]]_{i\pm\frac12} \}\}_{j\pm\frac12}}{\Delta x}+  \frac{ [[v]_{i\pm1}]_{j\pm1} }{\Delta y} \right ) \\
&\phantom{mmmmmmmm}=   \Delta x \, \del_x \mathscr D 
- \frac{ \Delta x^3}{12 } \del_x^3 \mathscr D 
+ \frac{ \Delta x^5}{120 } \del_x^5 \mathscr D 
- \frac{17  \Delta x^7}{20160 } \del_x^7 \mathscr D + \mathcal O(\Delta x^9)\\
&\phantom{mmmmmmmm}= 2 \frac{\exp(\Delta x \del_x) - 1}{\exp(\Delta x \del_x) + 1} \mathscr D
\end{align}
As the stationarity-consistency has been proven, one can be sure every order in this series to be a derivative of $\mathscr D$. Whenever $\mathscr D$ vanishes, the discrete second derivative vanishes as well. This calculation shows that the finite differences being ``rotationally-invariant'', as suggested in \cite{sidilkover02}, is not actually a relevant condition. They might seem so to first order in their expansion in powers of $\Delta x$, but they cannot remain so when higher order terms are taken into account.

The analysis easily generalizes to any number of spatial dimensions. These construction principles also allow one to derive stationarity preserving schemes of higher order.

\section{Conclusions and outlook}
\label{sec:conclusions}

This paper has presented a strategy for linear hyperbolic systems on how to study the existence of numerical stationary states for a given scheme (existence of a zero eigenvalue of the evolution matrix), and how these states can be found (their Fourier transform is the corresponding eigenvector).

The evolution of initial data that are a discretization of some analytically known stationary state therefore is to be seen as follows. Numerical stationary states are kept exactly constant in time, 
and all other solutions are stabilized by a certain amount of diffusion.
The difference between the initial data (e.g. derived from an analytic stationary state) and some representative of the numerical stationary states will be diffused away (or leave through the boundaries, if possible) and one will be left with the discrete stationary state. If the scheme is not stationarity preserving, then its stationary states discretize only a very small subset of all the analytic stationary states. This is, for example, the case with the upwind scheme. After some time basically all the features of the analytic stationary state will be diffused away and the resulting numerical stationary state will be unrecognizable. A stationarity preserving scheme on the other hand has stationary states that are discretizations of all the analytic stationary states. A comparison of simulation results with either type of schemes are shown in Fig. \ref{fig:greshostationary}--\ref{fig:greshostationaryfinal} and reflect the superiority of stationarity preserving schemes.

This paper presented examples and showed analysis of both dimensionally split and truly multidimensional stationarity preserving schemes for the acoustic equations.

The analysis of stationary states turns out to be a fruitful starting point in order to understand other properties of numerical schemes. For instance, stationarity preservation has been shown to be equivalent to the existence of a discrete vorticity that is kept exactly stationary during the numerical evolution. This allowed to make general statements about vorticity preserving schemes and to include schemes suggested previously in the literature into a larger framework. Additionally, it has been shown that these schemes are tightly linked (as described in sections \ref{sec:vorticity} and \ref{sec:lowmach}) to those that comply with the low Mach number limit for the acoustic equations. This connection has previously not been apparent.

Linearity is not a prerequisite for the existence of stationary solutions of hyperbolic PDEs. Subject of future work therefore is the development of stationarity-preserving methods for nonlinear equations. In the context of vorticity preservation, \cite{torrilhon04}, \cite{lerat07}, \cite{mishra09preprint} among others have already suggested examples of possible extensions. Future work intends to generalize the concepts of this paper. Additionally, with the tools developed in this paper it seems possible to study more complicated systems that contain source terms or have to fulfill constraints during the time evolution (structure preserving schemes).

\newcommand{\etalchar}[1]{$^{#1}$}

\end{document}